\documentclass[12pt]{amsart}
\usepackage{amsmath,amssymb}
\usepackage{mathpazo}
\usepackage{graphicx}

\usepackage{fullpage}
\def\RCS$#1: #2 ${\expandafter\def\csname RCS#1\endcsname{#2}}
\RCS$Revision: 1.28 $
\RCS$Date: 2005/05/13 20:15:46 $

\newcommand{\C}{{\mathbb C}}
\newcommand{\Q}{{\mathbb Q}}

\renewcommand{\H}{{\mathbb H}}
\newcommand{\Z}{{\mathbb Z}}
\newcommand{\Oint}{{\mathcal{O}}}
\newcommand{\maps}{\colon\thinspace}

\DeclareMathOperator{\Isom}{Isom}
\newcommand{\CP}{{ \mathbb{CP}}}
\newcommand{\PSL}[2]{\mathrm{PSL}_{#1} #2}
\newcommand{\SL}[2]{\mathrm{SL}_{#1} #2}
\newcommand{\GL}[2]{\mathrm{GL}_{#1} #2}
\newcommand{\spandef}[2]{{  \left\langle  {#1}  \ \left| \   {#2} \right. \right\rangle }}
\newcommand{\setdef}[2]{{  \left\{  \left. {#1}  \ \right| \   {#2} \right\} }}
\newcommand{\gpspan}[1]{{  \left\langle  {#1} \right\rangle }}
\newcommand{\mtext}[1]{\quad\mbox{#1}\quad}
\newcommand{\smallmat}[4]{{\left( \begin{array}{cc} {#1} & {#2} \\ {#3} & {#4} \end{array} \right)}}

\newcommand{\3}[1]{3\nobreakdash-\hspace{0pt}}

\makeatletter
\@ifundefined{swappedhead@plain}{%
\def\swappedhead#1#2#3{%
   \thmnumber{\@upn{\@secnumfont#2\@ifnotempty{#1}{.~}}}%
   \thmname{#1}%
    \thmnote{ {\the\thm@notefont(#3)}}}
\let\swappedhead@plain=\swappedhead}
\makeatother

\newtheoremstyle{plain}{}{}{\slshape}{}{\bfseries}{.}{0.5em}{}

\theoremstyle{plain} 

\swapnumbers
\newtheorem{theorem}{Theorem}[section]

\newtheorem{lemma}[theorem]{Lemma}

\newtheorem{requirement}[theorem]{Requirement}

\theoremstyle{definition}

\theoremstyle{remark}
\newtheorem{remark}[theorem]{Remark}

\makeatletter
  \let\c@theorem=\c@subsection
  \let\c@figure=\c@subsection
  \let\p@figure=\p@subsection
  
  \let\cl@figure=\cl@subsection
\makeatother

\begin{document}

\title{An ascending HNN extension of a free group inside $\SL{2}{\C}$}

\author[Calegari]{Danny Calegari} 
\address{Mathematics 253-37, California Institute of Technology, Pasadena, CA 91125, USA}
\email{dannyc@caltech.edu}

\author[Dunfield]{Nathan M.~Dunfield}
\address{Mathematics 253-37, California Institute of Technology, Pasadena, CA 91125, USA}
\email{dunfield@caltech.edu}

\begin{abstract} 
  We give an example of a subgroup of $\SL{2}{\C}$ which is a strictly
  ascending HNN extension of a non-abelian finitely generated free
  group $F$.  In particular, we exhibit a free group $F$ in
  $\SL{2}{\C}$ of rank $6$ which is conjugate to a proper subgroup of
  itself.  This answers positively a question of Drutu and Sapir
  \cite{DrutuSapir04}.  The main ingredient in our construction is a
  specific finite volume (noncompact) hyperbolic 3-manifold $M$ which
  is a surface bundle over the circle.  In particular, most of $F$
  comes from the fundamental group of a surface fiber.  A key feature
  of $M$ is that there is an element of $\pi_1(M)$ in $\SL{2}{\C}$ with
  an eigenvalue which is the square root of a rational integer.  We
  also use the Bass-Serre tree of a field with a discrete valuation to
  show that the group $F$ we construct is actually free.
\end{abstract}

\maketitle


\section{Introduction}
\label{sec:intro}

Suppose $\phi \maps F \to F$ is an injective homomorphism from a group
$F$ to itself.   The associated  HNN extension
\[
 H = \spandef{G, t}{ t g t^{-1} = \phi(g) \mtext{for} g \in G}
\]
is said to be \emph{ascending}; this extension is called
\emph{strictly ascending} if $\phi$ is not onto.  In
\cite{DrutuSapir04}, Drutu and Sapir give examples of residually finite
$1$-relator groups which are not linear; that is, they do not embed in
$\GL{n}{K}$ for any field $K$ and dimension $n$.  All of their examples are
ascending HNN extensions of free groups, and indeed this is how
they know the groups are residually finite (by \cite{BorisovSapir03}).
Motivated by looking at the $\SL{2}{\C}$ representations of their
examples, they asked whether $\SL{2}{\C}$ contains a strictly ascending
HNN extension of a non-abelian free group \cite[Problem
8]{DrutuSapir04}.  In particular, they asked if there is a non-abelian
free group $F$ in $\SL{2}{\C}$ which is conjugate to a proper
subgroup of itself.  In this note, we answer their question in the affirmative.
\begin{theorem}\label{main-thm}
  The group $\SL{2}{\C}$ contains a strictly ascending HNN extension
  of a free group of rank $6$.  
\end{theorem}
We first outline a general method of finding such groups, and then
describe a specific example in more detail.

\subsection{Surface bundles over $S^1$}   

A large part of our example comes from the fundamental group of a
hyperbolic \3-manifold which fibers over the circle.  Let $\Sigma$ be an
open surface of finite type, that is, a closed surface with a
non-empty finite set removed.  Assume $\Sigma$ has negative Euler
characteristic.  Given a homeomorphism $\phi:\Sigma \to \Sigma$ we can form the
mapping torus
\[
M_\phi = \raisebox{4pt}{$\Sigma \times [0,1]$} \Big/ \raisebox{-4pt}{$(s,1) \sim (\phi(s),0)$}
\]
Thurston showed that provided $\phi$ is pseudo-Anosov
(i.e.~$\phi$ is not homotopic to a finite order or reducible
homeomorphism), $M_\phi$ admits a {\em complete} hyperbolic
structure \cite{ThurstonFibered, Otal96}.  By Mostow rigidity, this
hyperbolic structure is in fact unique.

If $M_\phi$ is orientable, then after choosing an orientation, the
hyperbolic structure determines a discrete, faithful representation
\[
\rho:\pi_1(M_\phi) \to  \Isom^+(\H^3) = \PSL{2}{\C} .
\]
It turns out that $\rho$ always lifts to a representation into
$\SL{2}{\C}$ (see e.g. \cite{Culler86}), and henceforth we will regard
$\SL{2}{\C}$ as the target of $\rho$.

\subsection{The construction}

Now $\pi_1(M_\phi)$ is an HNN extension
\[
\pi_1(M_\phi) = \spandef{ \pi_1(\Sigma), t }{ t \alpha t^{-1} = \phi_*(\alpha) \text{ for all } \alpha \in \pi_1(\Sigma) }.
\]
Since $\Sigma$ is an open surface of finite type, its fundamental group
$\pi_1(\Sigma)$ is a finitely generated free group.  Take $F
= \rho(\pi_1(\Sigma))$ which is a free subgroup of $\SL{2}{\C}$.  Note
that for any $\beta \in \pi_1(M_\phi)$, the element $\rho(\beta)$ conjugates $F$
to itself.

We now restrict attention to $M_\phi$ with a very special property;
we will give an example of such an $M_\phi$ in
Section~\ref{sec:example}.
\begin{requirement}\label{requirement}
  There is a $\beta \in \pi_1(M_\phi) \setminus \pi_1(\Sigma)$ such that $\rho(\beta)$
  has trace $\pm (\sqrt{n} + 1/\sqrt{n})$ for some integer $n\geq 2$.
\end{requirement}
Equivalently, we want $\rho(\beta)$ to have eigenvalues $\pm \sqrt{n}$
and $\pm 1/\sqrt{n}$.  The image of $\rho(\beta)$ in $\SL{2}{\C}$ is
a hyperbolic isometry of $\H^3$, 
and has has two fixed points in $\partial \H^3 = \CP^1 = \C \cup \infty$.
After conjugating $\rho$, we can assume these fixed points are $0$ and
$\infty$, and $\rho(\beta)$ acts on $\CP^1$ by $z \mapsto nz$.  For $t \in \C$, let
$\mu_t$ denote the parabolic element of $\SL{2}{\C}$ which acts on
$\CP^1$ by $ z \mapsto z+t$.    Note that $\rho(\beta)$ conjugates $\mu_t$ to $(\mu_t)^n$

For each $t$, define $G_t = \gpspan{ F , \mu_t }$. Then $G_t$ is
finitely generated, and $\rho(\beta)$ conjugates it to a proper subgroup of
itself.  The following lemma, proved in the next section, shows that
$G_t$ is free for generic values of $t$.
\begin{lemma}\label{lemma}
  Let $F$ be a finitely generated subgroup of $\SL{2}{\C}$.  Consider
  the group $G_t = \langle F , \mu_t \rangle$ for some $t \in \C$.  Suppose that
  no non-trivial element of $F$ fixes $\infty \in \CP^1$.  Then for generic
  $t$, the group $G_t$ is isomorphic to the free product $F * \Z$.
\end{lemma}
For the action of the fundamental group of a finite volume hyperbolic
manifold on $\CP^1$, the stabilizer of any point is either empty, an
infinite cyclic group consisting of hyperbolic elements, or $\Z^2$
consisting of parabolic elements.  Now $\rho(\beta)$ is not parabolic as
its trace is not $\pm 2$.  Thus the stabilizer of $\infty$ in $\pi_1(M_\phi)$
is infinite cyclic.  So if some $f \in F$ also fixed $\infty$, then $f^n =
\rho(\beta)^m$ for suitable non-trivial powers $n$ and $m$.  Since no power
of $\beta$ lies in $\pi_1(\Sigma)$, the element $\rho(\beta)$ has no fixed point in
common with any element of $F$.  Thus the lemma applies in our context
and $G_t$ is free for most choices of $t$.

Thus, given a hyperbolic 3-manifold satisfying the requirement, we get
a free group $G_t$ which is conjugate to a proper subgroup of itself
by $\rho(\beta)$.  Then the subgroup $H = \gpspan{G_t, \rho(\beta)}$ of
$\SL{2}{\C}$ is isomorphic to the ascending HNN extension of $G_t$ by
the map induced via conjugation by $\rho(\beta)$.

This completes our construction and proves Theorem~\ref{main-thm}
modulo two things: proving the lemma, and exhibiting an $M_\phi$
satisfying Requirement~\ref{requirement} where $\pi_1(\Sigma)$ has rank 5.
These remaining tasks are carried out in the next two sections.  The
reader may wonder how plausible it is to expect a manifold satisfying
Requirement~\ref{requirement}.  For any finite volume hyperbolic
\3-manifold, the geometric representation $\rho \maps \pi_1(M) \to
\SL{2}{\C}$ can be conjugated so that its image lies in $\SL{2}{L}$
for some finite extension $L$ of $\Q$.  This is because if one looks
at representations $\pi_1(M_\phi) \to \SL{2}{\C}$ where the cusp group
acts by parabolics, then the representation $\rho$ coming from the
hyperbolic structure is an isolated point if we mod out by conjugation
\cite{GarlandRagunathan1970}.  The set of conjugacy classes of such
representations is an algebraic variety defined over $\Q$, and hence
isolated points have coordinates in a number field.  Thus it should
not be surprising that the trace of $\beta$ has a particular value in
a quadratic field.

\begin{remark}
  For our examples, the monomorphism of $G_t$ is necessarily
  reducible, i.e.~it preserves a free product decomposition of $G_t$
  into proper factors.  Mark Sapir asked whether there are also
  examples where this monomorphism is irreducible.
\end{remark}

\begin{remark}
  For our examples $G_t$ is necessarily indiscrete.  In fact this must
  always be the case.  Potyagailo and Ohshika proved that any
  topologically tame non-elementary 3-dimensional Kleinian group can
  not be conjugate in $\SL{2}{\C}$ to a proper subgroup
  \cite{OhshikaPotyagailo1998}.  Since all finitely generated Kleinian
  groups are topologically tame \cite{AgolTameness,
    CalegariGabaiTamness}, this implies that there are no examples
  where the free group base of the HNN extension is discrete.
\end{remark}

\subsection{Acknowledgments}  

Both authors were partially supported by the U.S.~National Science
Foundation (grant \#DMS-0405491) and the Sloan Foundation.

\section{Proof of the lemma}

\begin{proof}[{Proof of Lemma~\ref{lemma}}]
  The basic idea here is to use the Bass-Serre tree to show that the
  ``universal'' group $\gpspan{F, \mu_t}$, where $t$ is viewed as
  parameter, is isomorphic to the free product $F * \Z$.  It will then
  easily follow that $G_t$ is also isomorphic to $F * \Z$ for generic
  choices of $t$.
  
  For notational convenience, we change our parameterization of the
  parabolic $\mu_t$ from $z \mapsto z + t$ to $z \mapsto z + 1/t$.  In matrix form
  \[
  \mu_t =   \smallmat{ 1}{ t^{-1}}{ 0}{1}.
  \]
  Now let $K$ be the field of rational functions over $\C$.  Consider
  the homomorphism $i \maps F * \Z \to \SL{2}{K}$ which is just the
  inclusion of $\SL{2}{\C}$ into $\SL{2}{K}$ on the first factor, and
  sends the generator of $\Z$ to the above parabolic matrix.  We claim
  that $i$ is injective.
  
  Consider the discrete valuation on $K$ whose value on a rational function
  $f(t)$ is the order of vanishing of $f(t)$ at $t = 0$.
  Let $\Oint$ be the valuation ring of $K$; i.e.~the set of rational functions
  which do not have a pole at $0$.  Note that $\Oint$ is a local ring, and
  $t \Oint$ is its maximal ideal.  
  Then the action of $\SL{2}{K}$ on the Bass-Serre tree
  associated to this valuation gives us a splitting of $\SL{2}{K}$ as
  a free product with amalgamation as follows (see
  \cite[\S II.1.4]{SerreTrees} or \cite[\S3]{ShalenHandbook} for
  details).  Let $A = \SL{2}{\Oint}$ and $B = X A X^{-1}$ be the conjugate
  of $A$ by
  \[
  X = \smallmat{t^{-1}}  0  0 1  \in \GL{2}{K}, \text{ i.e. } B = \smallmat  a {t^{-1}b} { t c }  d  \text{ where } \smallmat a b c d \in \SL{2}{\Oint}.
  \]
  If we set
  \[
  C = \setdef{ \smallmat a b c d \in \SL{2}{\Oint} }{  c \in t \Oint},
  \]
  then $\SL{2}{K} = A *_C B$.   (See the top of page 81 of
  \cite[\S II.1.4]{SerreTrees} for this precise splitting in the analogous case of
  $\SL{2}{\Z[1/p]}$.)  
  
  Note that the map $i \maps F * \Z \to \SL{2}{K}$ takes $F$ into the
  first factor of the amalgam and takes $\Z$ into the other.  The
  condition that no element of $F$ fix $\infty$ exactly corresponds to the
  condition that the image $i(F)$ in $\SL{2}{K}$ does not intersect
  the edge group $C$ joining the two factors of the amalgam.  As
  $i(\Z)$ is also disjoint from $C$, the subgroup of $\SL{2}{K}$
  generated by $i(F)$ and $i(\Z)$ is isomorphic to $i(F) * i(\Z)$.
  Thus $i$ is an isomorphism.
  
  To complete the proof of the lemma, we will show that for a generic
  choice of $t$, the composite of $i$ with the induced map $\SL{2}{K}
  \to \SL{2}{\C}$ is also injective.  Let $w \in F * \Z$ and consider
  the set of $t$ for which $w$ maps to the identity under the map $F *
  \Z \to G_t$.  Such $t$ are the solutions to some polynomial equations
  in $t$; as $i$ is injective, these equations are nontrivial and have
  a finite solution set.  As $F$ is finitely generated, there are only
  countably many such $w \in F * \Z$ and hence only countably many
  ``bad'' values for $t$.  Thus if we select $t \in \C$ from outside
  this countable set we have that $G_t$ is the free product $F * \Z$,
  as desired.
\end{proof}

\begin{remark}  
  If, as is the case in our application, the entries of the elements
  of $F$ live in some number field $L$, then one can simply take $t$
  to be any transcendental number, e.g. $\pi$ or $e$.
\end{remark}

\section{Example}\label{sec:example}

Our example of a fibered $3$-manifold satisfying
Requirement~\ref{requirement} is the orientable double cover of the
hyperbolic manifold $M = y505$ from the Callahan-Hildebrand-Weeks
census \cite{CallahanHildebrandWeeks} that comes with SnapPea
\cite{SnapPea}.  The manifold $M$ is nonorientable and has one torus
cusp.  It has an ideal triangulation with $7$ tetrahedra.  The shapes of
these tetrahedra all lie in the field $\Q(e^{\pi i/3})$.  The
fundamental group of $M$ has a presentation
\[
\pi_1(M) = \spandef{ a,b} {aabbaaBAbaBabAABBAbAAB = 1}
\]
where $A = a^{-1}$ and $B = b^{-1}$.  The \emph{faithful} action of
$\pi_1(M)$ on $\CP^1$ is given by:
\[
a \maps z \mapsto   1 + \frac{1 + \omega^2}{ \overline{z}} \text{ and } b \maps z \mapsto  \frac{  - 3 \overline{\omega} z + 9 \omega - 6 \overline\omega}{  (\omega + \overline\omega) z - 7 \omega + 5 \overline\omega}
\]
where $\omega = e^{\pi i/6}$, a primitive $12^{\mathrm{th}}$ root of unity.
Notice here that $a$ is orientation reversing, and that expressed as a matrix
\[
a^2 = \frac{1}{\sqrt{3}} \smallmat{2 + \omega^2}{2 - \omega^2}{1}{2 - \omega^2}
\]
which has trace $\sqrt{3} + 1/\sqrt{3}$, as desired.

In a minute, we will pass to the orientation cover $N$ of $M$ to
produce the needed example, but first let's observe that $M$ fibers
over the circle.  First notice that the defining relation for
$\pi_1(M)$ is in the commutator subgroup of the free group on $a$ and
$b$, and hence $H_1(M; \Z) = \Z \oplus \Z$.  Consider the homomorphism $\phi
\maps \pi_1(M) \to \Z$ given by $\phi(a) = 1$ and $\phi(b) = 0$.  If we
set $b_k = a^k b a^{-k}$ and $B_k = b^{-1}_k$ we see that the defining
relation of $\pi_1(M)$ becomes
\[
b^2_2 B_4 b_3 B_4 b_5 B_3^2 b_2 B_0 = 1.
\]
Since the highest and lowest subscripts appear only once each, we see
by Magnus rewriting \cite{MagnusKarrassSolitar} that the kernel of
$\phi$ is freely generated by $\{ b_0, b_1, b_2, b_3, b_4 \}$.  Thus by
a theorem of Stallings \cite{Stallings62}, the map $M \to S^1$ induced
by $\phi$ can be homotoped to a fibration where the fundamental group of
the fiber is the kernel of $\phi$.  (Actually, for
Theorem~\ref{main-thm}, we don't really need Stallings' theorem since
we have explicitly exhibited the desired free normal subgroup.
However, ultimately the justification that this subgroup of matrices
is free comes from identifying the full group of matrices with
$\pi_1(M)$; that is, the representation given above is faithful.  This
is seen by deriving the representation from a hyperbolic structure on
$M$, described by a geometric solution to the gluing equations for
some ideal triangulation.)
\begin{figure}
\begin{center}
  \includegraphics[scale=0.8]{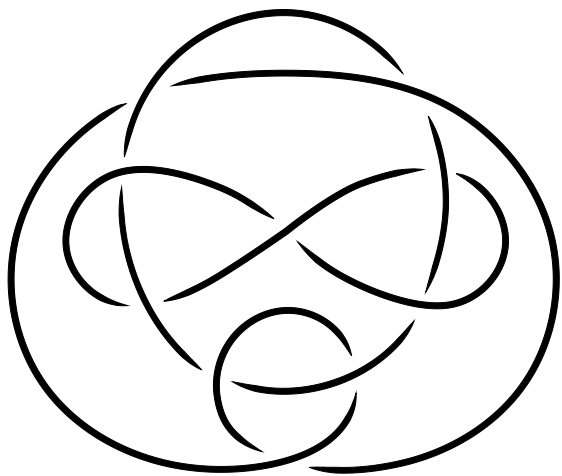}
\end{center}
\caption{Doing $0$ surgery on the circular component with no self-crossings yields the manifold $N$.   This is link $10^3_5$ in the Christy table of links through $10$ crossings.}\label{fig}
\end{figure}

Now pass to the orientable double cover $N$ of $M$, whose fundamental
group is $\phi^{-1} ( 2 \Z)$.  (For a Dehn surgery description of $N$,
see Figure~\ref{fig}.)  Note that the kernel of $\phi$ restricted to
$\pi_1(N)$ hasn't changed and is still the free group $F$ of rank 5.
Moreover, $a^2$ is not in $F$ as $\phi(a^2) = 2$ and, as we mentioned,
the trace of $a^2$ is $\sqrt{3} + 1/\sqrt{3}$.  Thus $N$ satisfies
Requirement~\ref{requirement}, and this completes the proof of
Theorem~\ref{main-thm}.

\begin{remark}
  It is easy to use Brown's elegant algorithm \cite{Brown87} for
  computing the Bieri-Neumann-Strebel invariant of a 1-relator group
  to see that any $\phi \maps \pi_1(M) \to \Z$ except $b^*$ and $a^* - 2
  b^*$ corresponds to a fibration.  (Here $\{ a^*, b^*\}$ is the dual
  basis of $H^1(M; \Z)$ with respect to the basis $\{a, b\}$ of
  $H_1(M; \Z)$.)  This gives infinitely many other examples, where the
  rank of the free group in question grows more slowly than in examples
  obtained by taking subgroups of finite index in the kernel
  of $a^*$.
\end{remark}

\begin{remark}
  Of the $6070$ cusped hyperbolic 3-manifolds in the
  Callahan-Hildebrand-Weeks census, $y505$ seems to be the only one
  that satisfies Requirement~\ref{requirement}.  There are a handful of
  other examples (mostly also non-orientable) which have an element
  with the right trace, but these are not fibered. However, it is
  plausible that they have finite covers which are fibered.
\end{remark}

\providecommand{\bysame}{\leavevmode\hbox to3em{\hrulefill}\thinspace}
\providecommand{\MR}{\relax\ifhmode\unskip\space\fi MR }
\providecommand{\MRhref}[2]{%
  \href{http://www.ams.org/mathscinet-getitem?mr=#1}{#2}
}
\providecommand{\href}[2]{#2}

\end{document}